\documentstyle[eqsection]{article}

\newenvironment{pf}{\noindent{\sc Proof}.\enspace}{\rule{2mm}{2mm}\medskip}

\newtheorem{theorem}{Theorem}[section]
\newtheorem{lemma}{Lemma}[section]

\newtheorem{condition}{Condition}[section]
\newtheorem{remark}{Remark}[section]
\newtheorem{remarks}{Remark}[section]
\newtheorem{definition}{Definition}[section]
\newcommand{\be}{\begin{equation}}
\newcommand{\ee}{\end{equation}}

\newcommand{\teta}{\theta}

\newcommand{\om}{\omega}

\newcommand{\ep}{\epsilon}

\newcommand{\ov}{\overline}

\newcommand{\wtilde}{\widetilde}

\renewcommand{\a }{\alpha }
\renewcommand{\b }{\beta }
\newcommand{\s }{\sigma }
\renewcommand{\d }{\delta }

\newcommand{\e }{\varepsilon }
\newcommand{\g }{\gamma}

\newcommand{\vphi}{\varphi }

\renewcommand{\t }{\tau }

\newcommand{\dps}{\displaystyle}

\begin{document}

\title{{\bf Fast Arnold Diffusion in three time scale systems}}

\author{Massimiliano Berti and Philippe Bolle}
\date{}
\maketitle

{\bf Abstract:}
We consider the problem of Arnold Diffusion for 
nearly integrable partially isochronous Hamiltonian 
systems with three time scales.  
By means of a careful shadowing analysis, based on a variational technique,
we prove that, along special directions, Arnold diffusion 
takes place with fast (polynomial) speed, even though the ``splitting 
determinant'' is exponentially small.
\footnote{Supported by M.U.R.S.T. Variational Methods and Nonlinear
Differential Equations.}
\\[2mm]

\section{Introduction}
In a previous paper \cite{BB3} (see also \cite{BBN}) we introduced,
in the context of nearly integrable Hamiltonian systems,
a functional analysis approach to the ``splitting 
of separatrices'' and to the ``shadowing problem''.
We applied our method to the problem 
of Arnold Diffusion, i.e. topological instability 
of action variables, for nearly integrable partially 
isochronous systems.
The aim of this paper is to improve the shadowing theorem of
\cite{BB3}
and to apply this new theorem 
to the three time scale system (\ref{eq:Hamil}) below, in order
to prove that 
along special directions Arnold diffusion takes 
place with ``very fast speed'',
namely a speed  polynomial in $ \e $. To that effect,
we  use the results on the splitting provided in \cite{BB3}.
 
Three time scale Hamiltonian systems have been introduced 
in \cite{CG} as a description of the D'Alembert problem in 
Celestial Mechanics. Later on three time scale systems
have been reconsidered for example in    
\cite{GGM}, \cite{GGM3}, \cite{PV}, \cite{BCV}, \cite{BB3}, 
\cite{GGMt}.

In this paper we focus on isochronous three time scale systems as
\be\label{eq:Hamil}
{\cal H}_\mu =  \frac{1}{\sqrt{\e}} I_1 + \e^a \b \cdot I_2 
+ \frac{p^2}{2} + ( \cos q - 1) ( 1 + \mu f( \vphi )),  
\ee
where $ (\vphi_1, \vphi_2, q) \in 
{\bf T}^1 \times {\bf T}^{n-1} \times {\bf T}^1 $ are 
the angle variables,  $ ( I_1, I_2, p) 
\in {\bf R}^1 \times {\bf R}^{n-1} \times {\bf R}^1 $
are the action variables, $ \b =(\b_2, \ldots, \b_n) \in {\bf R}^{n-1} $,
$ n \geq 3 $, $ a  > 0 $ and 
$ \e > 0 $, $ \mu \geq 0 $ are small real parameters.
We will assume that $ \mu = O( \min \{  \e^{ 3/2}, \e^{2 a + 1} \})$. 
Hamiltonian $ {\cal H}_\mu $ describes a
system of $ n $ isochronous harmonic oscillators
with a Diophantine frequency vector $ \om_\e = (1 / \sqrt{\e}, \e^a \beta )$ ,
with one fast frequency $ \om_{\e, 1} = 1/ \sqrt{\e}$ and 
$(n - 1)$ slow frequencies $ \om_{\e, 2} =  \e^a \b $,  
weakly coupled with a pendulum.
\\[1mm]
\indent
When $ \mu = 0 $ the energy $ \om_{\e, i} I_i $ of each oscillator is a
constant of the  motion.
The problem of {\it Arnold diffusion} in this context is whether,
for $ \mu \neq 0 $, there exist motions whose net effect is to transfer
$O(1)$-energy from one oscillator to others in a certain time $ T_d $ called 
the diffusion time. 
\\[1mm]
\indent
The existence of Arnold diffusion is usually proved 
following the mechanism proposed in \cite{Arn}.
For $ \mu = 0 $ Hamiltonian ${\cal H}_\mu$ admits
a continuous family of $ n $-dimensional partially hyperbolic invariant tori
${\cal T}_{I_0} = \{ (\vphi , I, q,p) \in  {\bf T}^n
\times  {\bf R}^n \times {\bf T}^1 \times {\bf R}^1 \ | 
\  I = I_0, \ q = p = 0 \}$
possessing stable and unstable manifolds
$ W^s ({\cal T}_{I_0}) = W^u({\cal T}_{I_0}) = 
\{ (\vphi , I, q,p) \in  {\bf T}^n
\times  {\bf R}^n \times {\bf T}^1 \times {\bf R}^1 \  |
  \ I = I_0, \ p^2 / 2 + ( \cos q - 1) = 0 \}$
called ``whiskers'' by Arnold.
For $ \mu  $ small enough the perturbed stable and unstable manifolds
$ W^s_\mu ({\cal T}_{I_0}^\mu ) $ and
$ W^u_\mu ({\cal T}_{I_0}^\mu ) $ may split and intersect transversally,
giving rise to a chain of tori connected by heteroclinic orbits.
By a shadowing type argument one can then prove the existence of an orbit
such that the action variables $ I $ undergo a variation of $ O (1) $
in a certain time $ T_d $ called the {\it diffusion time}.
In order to prove the existence of diffusion orbits
following the previous mechanism one encounters two different problems:
$1)$ {\it Splitting of the whiskers}; $2)$ {\it Shadowing problem}.

The ``splitting of the whiskers'' for Hamiltonian $ {\cal H}_\mu $, 
when  $ \mu = O( \e^p )$, $ p > 0 $ and $ \e \to 0 $, 
has been studied in \cite{GGM}, \cite{GGMt}, \cite{PV} and \cite{BB3}.
In \cite{GGM}-\cite{GGMt} and \cite{PV} the size of the splitting is measured
by the ``determinant of the splitting matrix'' which 
turns out to be exponentially small, precisely  $ O( \exp (- (\pi /2) \e^{-1/2})) $.
We underline that papers \cite{GGM}-\cite{GGMt} deal also
with non-isochronous systems and more general perturbation terms
(but two rotators only).

In \cite{BB3},
the splitting of stable and unstable manifolds 
is related to  the variations of  the ``homoclinic
function'' $ G_\mu: {\bf T}^n \to {\bf R} $ (defined in (\ref{eq:homo})),
which is the difference between the generating functions
of stable and unstable manifolds at section $ \{ q = \pi \} $. 
$ \nabla G_\mu (A) $ provides a measure of  the distance between 
stable and unstable manifolds, so that 
a critical point $ \ov{A} $ of $ G_\mu $ gives rise to a homoclinic 
intersection. Usually  det $ D^2 G_\mu ( \ov {A} ) $ is called
the ``splitting determinant''.
The use of the ``homoclinic function'' $ G_\mu $ for measuring the
splitting has two advantages. Firstly, it
is very well suited to deal with the shadowing problem  by means
of variational  techniques
because $ G_\mu $  is nothing but the difference of the values
of  the Lagrangian action functional
associated to  the quasi-periodically forced pendulum (\ref{pendper}) 
at two true  solutions, lying respectively on the stable and unstable manifolds
$ W^{s,u}_\mu ( {\cal T}_{I_0} )$, see (\ref{def:homo}).
Secondly 
it may shed light on a ``non uniform'' splitting which would  not be
given by the splitting determinant, when the 
variations of $ G_\mu $ in different directions are of different orders. 

For the three time scale system associated to Hamiltonian
${\cal H}_\mu$, ``non uniform'' splitting is suggested by
the behaviour of the first order expansion of $ G_\mu $ in $\mu$, called the 
Poincar\'e-Melnikov approximation.  
In fact the first order term, which is  given by the 
Poincar\'e-Melnikov primitive defined in (\ref{PoiMel}), has 
exponentially small oscillations in the fast angle $ A_1 $,
and polynomially small ones in the slow angles $ A_2 $.
Naively this hints the splitting
to be exponentially small in the direction
$ I_1 $ and just polynomially small 
in the directions $ I_2 $.

However, in general, for $ \mu = O ( \e^p ) $ and $ \e \to 0 $
the homoclinic function $ G_\mu $ 
is not well approximated by the Poincar\'e-Melnikov primitive.
In \cite{GGM}-\cite{GGMt} the asymptotic validity of Melnikov's integrals 
for computing the exponentially small ``splitting determinant'' is 
proved to hold only after exhibiting many cancellations.

In \cite{BB3} the naive Poincar\'e-Melnikov approximation
for Hamiltonian ${\cal H}_\mu $ has been rigourously justified 
for $\mu \e^{- 3/2} $ sufficiently small, in a different way. 
We define another ``splitting function''
$ \wtilde{G}_\mu $, see (\ref{def:wtildeg}),
whose critical points as well give rise to homoclinic 
intersections.  $ \wtilde{G}_\mu $ 
is well approximated, for $\mu = O( \e^p )$ and $ \e \to 0 $,  
by the Poincar\'e-Melnikov primitive and has 
exponentially small oscillations in $ A_1 $, see theorem
\ref{thm:tts}.
The crucial observation is that 
$ G_\mu $ and $ \wtilde{G}_\mu $ are the same function
up to a diffeomorphism $ \psi_\mu $ of the torus close to identity, namely 
$ \wtilde{G}_\mu = {G}_\mu \circ \psi_\mu  $, see theorem 
\ref{thm:diffeo}. 
\\[2mm]
\indent
After the works \cite{Bs}, \cite{Bs1}, \cite{M}, \cite{Cr0},  
\cite{BCV}, \cite{BB3}, \cite{Cr}, \cite{CrG} and references therein,  
it is a well established fact that the diffusion time 
is estimated by a polynomial inverse power of the splitting.
For instance,
using the estimate on the size of the splitting of \cite{GGM} and
\cite{GGMt} an exponentially long diffusion time has been obtained in
\cite{BCV}, namely 
$ T_d = O(\exp (C / \e^{b})) $ for some $ b > 0 $   
(see also theorem 5.2 of \cite{BB3}). 

However  the properties of $G_\mu$ (oscillations of different
amplitude orders according to the  direction) suggest that 
Arnold diffusion can take place with different speed along different 
directions; since, for larger splitting one would expect a
faster speed of diffusion, one could guess the existence of diffusion 
orbits that drift along the ``fast'' directions 
$ I_2 \in {\bf R}^{n-1} $,   
where the splitting is just polynomially small w.r.t. $1/ \e$,
in a polynomially long diffusion time $ T_d = O( 1/ \e^q )$.
The aim of this paper is to prove that this is indeed the case.
In order to prove this phenomenon (see {\it theorem \ref{step1}} for the 
general case and {\it theorem \ref{thm:main}} for an application) 
we refine the shadowing theorem 2.3 of \cite{BB3} 
for dealing with the present ``non-uniform'' splitting.
Note that, because of the preservation of the energy along
the orbits, Arnold diffusion can take place in the direction
$I_2$ for $n\geq 3$ only.
\\[1mm]
\indent
In order to justify heuristically our result
we recall how the diffusion time $ T_d $ is estimated in \cite{BB3},
once it is verified that stable and unstable manifolds split.
$ T_d $ is, roughly, estimated by the product
of the number of heteroclinic transitions $ k $ ($=$ number of tori
forming the transition chain = heteroclinic jump/splitting)
and of the time $ T_s $ required for a single transition,
namely $ T_d = k T_s $. 
The time for a single transition $ T_s $  is bounded by the maximum
time between the   ``ergodization time'' $T_e$ of the torus 
${\bf T}^n $ run by the linear flow $ \om_\e t $,
and the time needed  to ``shadow''
homoclinic orbits for the corresponding quasi-periodically forced pendulum
equation \ref{pendper}.

The reasons for which we are able to move in polynomial time w.r.t $ 1/\e $ 
along the fast $ I_2 $ directions are the following three ones. 
($i$) As in \cite{BB3}, since the homoclinic orbit 
decays exponentially fast to $0$, the time needed  to ``shadow''
homoclinic orbits for the quasi-periodically forced pendulum (\ref{pendper})
is only polynomial.  ($ii$)
Since the splitting is polynomially small 
in the directions $ I_2 $, 
 we can choose just a polynomially large number of tori 
forming the transition chain $ k = O(1/ \e^p)$ to get 
a $O(1)$-drift of $I_2$. ($iii$)
Finally, the most difficult task is getting a  polynomial estimate for the 
``ergodization time'' $T_e$ 
-defined as the time needed for the flow 
$ \{ \om_\e t \}$ to make an $ \alpha $-net of the torus- with
$\alpha$ appropriately small. 
By a result of \cite{BGW} this time satisfies 
$ T_e = O ( 1 / \alpha^\tau ) $.
Let us explain  how this estimate enters into play.
In order to apply our ``gluing'' variational technique, 
the projection of our shadowing orbit on the torus $ { \bf T }^n $,
namely $ \{ \om_\e t + A_0 \} $, must approach, at each transition,
sufficiently close to the homoclinic point $ \ov{A} $ to be capable
to ``see'' the homoclinic critical point $ \ov{A} $ of $ G_\mu $. 
The crucial improvement of the shadowing theorem \ref{step1}
allows the shadowing orbit to approach $\ov{A}$ 
only up to a polynomially small distance $ \alpha = O( \e^p ) $, $ p >0 $, 
(and not exponentially small as it 
would be required when applying the shadowing theorem of \cite{BB3}).
By the forementioned estimate on the ergodization time 
$ T_e = O ( 1 / \alpha^\tau ) $ it results that the 
minimum time after which the homoclinic trajectory 
can ``jump'' to another torus is only polynomially long w.r.t $1/ \e$.
Actually this allows to improve as well the exponential estimate 
on the diffusion time required to move also in the $I_1$ 
direction, see remark \ref{poliexp}.
\\[1mm]
Theorems \ref{step1} and \ref{thm:main} are 
the first step to prove the existence of this phenomenon also
for more general systems (with non isochronous terms 
and more general perturbations).

The paper is organized as follows:
in section 2 we recall some preliminary results taken from \cite{BB3}.  
In section 3 we introduce the general ``splitting condition'' 
which will be used in section 4 to prove the shadowing theorems.
\\[2mm]
Through the paper $ C_i $ and $ \d_i $ will denote positive constants which 
are independent of $ \e $ and $\mu$. 

\section{Preliminaries}

In this section we recall the results of \cite{BB3}  
that will be used in the sequel. 
We refer to \cite{BB3} for complete details and for 
the description of the general functional analysis approach  
based on a Lyapunov Schmidt type reduction. With respect to the
notations of \cite{BB3} we remark that we 
have changed the sign of the perturbation $f$
in Hamiltonian $ {\cal H}_\mu $.
\\[1mm]
\indent
The equations of motion derived by Hamiltonian $ {\cal H}_\mu $ are
\be\label{eqmotion}
\dot{\vphi} = \om_\e, \qquad
\dot{I} = - \mu (\cos q -1 ) \ \partial_{\vphi}f(\vphi ), \qquad
\dot{q} = p, \qquad
\dot{p} = \sin{q} \ (1 + \mu \ f (\vphi )).
\ee
The angles $ \vphi $ evolve as $ \vphi (t) =   \om_\e t + A $;
 therefore  equations (\ref{eqmotion})
can be reduced to the quasi-periodically forced pendulum
equation
\be\label{pendper}
- \ddot{q} + \sin{q} \ (1 + \mu f ( \om_\e t + A)) = 0,
\ee
corresponding to the Lagrangian
\be\label{lagraper}
{\cal L}_{\mu, A} (q, \dot{q},t) = \frac{{\dot q}^2}{2} + 
(1- \cos q) (1 + \mu f(\om_\e t + A)).
\ee
For each solution $ q(t) $ of (\ref{pendper}) one recovers the dynamics
of the actions $ I (t) $ by quadratures in (\ref{eqmotion}).
\\[1mm]
\indent
For $ \mu = 0 $ equation (\ref{pendper})
possesses the one parameter family 
of homoclinic solutions to $ 0 $, mod $ 2 \pi $,
$ q_\teta (t) =$ $4 \ {\rm arctan} (\exp {(t - \teta)}),
\  \teta \in {\bf R} $.
Using the Implicit Function Theorem 
one can prove (lemma 2.1 of \cite{BB3}) that 
there exist, near the unperturbed homoclinic solutions $ q_\teta (t) $, 
for $ 0< \mu < \mu_0 $ small enough independently of $ \om_\e $, 
{\it ``pseudo-homoclinic solutions''}
$ q_{A,\teta}^\mu (t)$ of equation (\ref{pendper}).
These are true solutions of (\ref{pendper})
in each interval $(- \infty, \teta)$ and $(\teta, +\infty)$;
at time $ t = \teta $ such pseudo-solutions are glued with continuity
at value  $ q^\mu_{A, \teta} ( \teta ) = \pi $ and for
$ t \to \pm \infty $ are asymptotic to the equilibrium $0 $ mod $ 2 \pi $.
We can then define the function
$ F_\mu : {\bf T}^n \times {\bf R} \to {\bf R} $ as the action 
functional of Lagrangian (\ref{lagraper})
evaluated on the ``1-bump pseudo-homoclinic solutions''
$ q_{A, \teta }^\mu (t) $, namely
\be\label{def:homo}
F_\mu (A, \teta) =
\int_{ - \infty}^\teta {\cal L}_{\mu, A} ({  q}_{A,\teta}^\mu (t),
\dot{  q}^\mu_{A,\teta} (t), t) \ dt
+ \int_\teta^{+ \infty}
{\cal L}_{\mu, A} ({  q}_{A,\teta}^\mu (t),
\dot{  q}^\mu_{A, \teta } (t),t) \ dt,
\ee
and the {\it ``homoclinic function''} $ G_\mu : {\bf T}^n \to {\bf R} $ as
\be\label{eq:homo}
G_\mu ( A ) = F_\mu ( A, 0 ).
\ee
There holds 
\be \label{invpro}  F_\mu  (A , \teta ) = G_\mu (A + \om_\e \teta), \forall
\teta \in {\bf R}. \ee

\begin{remark}
The homoclinic function $ G_\mu $ is the difference
between the generating functions
$ {\cal S}_{\mu,I_0}^\pm (A, q )$ of the stable and the unstable manifolds
$ W_\mu^{s,u} ( {\cal T}_{I_0} )$
(which in this case are {\it exact} Lagrangian manifolds)
at the {\it fixed}  section $ \{ q = \pi \} $, namely
$G_\mu (A) = {\cal S}_{\mu, I_0}^- (A, \pi ) -
{\cal S}_{\mu, I_0}^+ ( A, \pi ) $. 
A critical point of $ G_\mu $ gives rise to a homoclinic orbit
to torus ${\cal T}_{I_0}$, see lemma 2.3 of \cite{BB3}. 
\end{remark}

In order to justify the dominance
of the Poincar\'e-Melnikov function when $ \mu = O(\e^p )$
one would need to extend analytically the function $F_\mu (A, \teta)$ for
complex values of the variables.
Since the condition $ q^\mu_{A, \teta } ( Re \ \teta ) = \pi $, appearing
naturally when trying to extend the definition of 
$ q^\mu_{A, \teta}$ to $\theta \in {\bf C}$, breaks
analyticity, the function $ F_\mu (A, \teta) $ can not 
be easily analytically extended in a 
sufficiently wide complex strip. To overcome this problem, in \cite{BB3}
the Lagrangian action functional is 
evaluated on different ``1-bump pseudo-homoclinic solutions''
$ Q^\mu_{A,\teta} $. 
Define $ \psi_0: {\bf R} \to {\bf R} $ by
$ \psi_0 (t) = \cosh^2 ( t ) / ( 1 + \cosh t )^3 $ and
set $ \psi_\teta (t)= \psi_0 (t - \teta)$.
Two important properties of the function $ \psi_0 (t) $ are that 
$ \int_{\bf R} \psi_0 (t) {\dot q}_0 (t) \ dt \neq 0 $ and that 
it can be extended
to a holomorphic function on ${\bf R} + i( - \pi , \pi) $
(while the homoclinic solution $ q_0 (t) $ can be extended to 
a holomorphic function only up to ${\bf R} + i( - \pi / 2, \pi / 2) $).
By the Contraction Mapping Theorem there exist (lemma 4.1 of \cite{BB3})
near $ q_\teta $, for $ \mu $ small enough,
{\it pseudo-homoclinic solutions} $ Q_{A, \teta}^\mu (t) $
and a constant $ \a_{A, \teta}^\mu $ defined by
$$
- \ddot{Q}_{A, \teta}^\mu + \sin{Q}_{A, \teta}^\mu ( 1 + \mu \
f (\om_\e t + A) ) = \a_{A, \teta}^\mu \psi_\teta (t) 
\quad {\rm and} \quad
\int_{\bf R} \Big( Q_{A, \teta}^\mu (t)- q_\teta (t) \Big)
\psi_\teta (t) \ dt = 0.
$$
We define the function
$ \wtilde{F}_\mu: {\bf T}^n \times {\bf R} \to {\bf R} $
as the action functional of Lagrangian  (\ref{lagraper})
evaluated on the ``1-bump pseudo-homoclinic solutions''
$ Q_{A, \teta }^\mu (t) $, namely
\be\label{def:wtildeg}
\wtilde{F}_\mu (A, \teta) =
\int_{\bf R}  {\cal L}_{\mu, A} ( Q_{A,\teta}^\mu (t),
\dot Q^\mu_{A,\teta} (t), t) \ dt
\ee
and $ \wtilde{G}_\mu: {\bf T}^n  \to {\bf R} $ as
$ \wtilde{G}_\mu ( A ) = \wtilde{F}_\mu ( A, 0 ).$

\begin{remark}
Also critical points of $ \wtilde{G}_\mu $
give rise to homoclinic solutions to torus ${\cal T}_{I_0}$, 
see lemma 4.2 of \cite{BB3}.
By theorem \ref{thm:diffeo} below, from a geometrical point of view the introduction of the ``homoclinic
function'' ${\wtilde G}_\mu $ may be interpreted simply as measuring the
splitting with a non constant Poincar\'e section, see 
the introduction of \cite{BB3}.
\end{remark}

The crucial point is now to observe that 
the homoclinic functions  $ G_\mu $ and $ \wtilde{G}_\mu $  
are the same up to a change of variables close to the identity,
as stated by the following theorem (see theorem 4.1 of \cite{BB3})

\begin{theorem}\label{thm:diffeo}
For $ \mu $ small enough  (independently of $ \om_\e $) 
there exists a Lipschitz homeomorphism (a real analytic diffeomorphism 
if $ f $ is analytic) $ \psi_\mu : {\bf T}^n \to {\bf T}^n $ of the form
$ \psi_\mu (A) = A + k_\mu (A) \om_\e $ with
$ k_\mu: {\bf T}^n \to {\bf R} $ satisfying
$ k_\mu ( A ) = O ( \mu ) $, $ | k_\mu ( A ) - k_\mu ( A' )| = 
O(\mu ) | A - A' | $ 
such that $  {\wtilde G}_\mu = G_\mu \circ \psi_\mu.$
\end{theorem}

Let $\Gamma (\e, A)$ denote the Poincar\'e-Melnikov primitive
\be\label{PoiMel}
\Gamma (\e, A) =
\int_{{\bf R}} (1- \cos q_0 (t)) f (\om_\e t + A) \ dt.
\ee
Develop in Fourier series w.r.t. the first variable 
the  homoclinic function 
$ \wtilde{G}_\mu (A) = \sum_{k_1 \in {\bf Z}}
\wtilde{g}_{k_1} (A_2) e^{i k_1 \cdot A_1 }$ 
and the Poincar\'e-Melnikov primitive
$ \Gamma (\e, A) = \sum_{k_1 \in {\bf Z}}
\Gamma_{k_1} (\e, A_2 ) e^{i k_1 \cdot A_1  }$.
Assume that the perturbation $f$ is analytic w.r.t 
$( \vphi_2 , \ldots, \vphi_n )$. More precisely assume that
there exist $ r_i > 0 $ for $ i = 2, \dots, n $,  such that
$ f $ has a $ C^{\infty} $ extension in
$ D:=  {\bf R} \times ({\bf R} + i [ -r_2, r_2 ]) \times \ldots \times
({\bf R} + i [ -r_n, r_n ])$,
holomorphic w.r.t. $( \vphi_2 , \ldots, \vphi_n )$.
Denote the supremum of $ | f | $ over $ D $ as  
$ || f || := \sup_{ \vphi \in D} | f ( \vphi )|.$
The following theorem about the splitting of stable 
and unstable manifolds in three time scale systems, holds
(see theorem 5.1 of \cite{BB3})

\begin{theorem}\label{thm:tts}
For $ \mu ||f|| \e^{- 3/2} $ small there holds, for all 
$ A_2 \in {\bf T}^{ n - 1} $ 
\begin{eqnarray*}
\wtilde{G}_\mu (A_1, A_2) & = & 
\wtilde{g}_0 ( A_2 ) + 2  {\rm Re} \Big[ \wtilde{g}_1 ( A_2 ) e^{i A_1} \Big]
+ \wtilde{R} (A_1, A_2 )\\
& = & Const + 
\Big( \mu \Gamma_0 (\e, A_2 ) + R_0 (\e, \mu, A_2 ) \Big) 
+  2 \ {\rm Re} \ 
\Big[ \mu \Gamma_1 (\e, A _2) + R_1 (\e, \mu, A_2) \Big] e^{i A_1} \\
& + & \wtilde{R} (A_1, A_2 ), 
\end{eqnarray*}
where
$ R_0 (\e, \mu,  A_2) = O \Big( \mu^2 ||f||^2 \Big), \quad 
R_1 (\e, \mu, A_2) = O \Big( \dps \frac{ \mu^2||f||^2 }{\e^2}  
\exp{ \Big( -\dps \frac{\pi}{2 \sqrt{\e}} \Big)} \Big),$ and 
$$
\wtilde{R} (A_1, A_2 ) =
O \Big( \mu  \e^{-1/2} ||f|| \exp{\Big(- \frac{\pi}{\sqrt{\e}} \Big)} \Big).
$$
\end{theorem}

In order to prove our shadowing theorem we need also to recall 
the definition of the ${\bf k}${\bf -bump pseudo-homoclinic solutions}
$ q_{A,\teta}^L (t)$ 
for the quasi-periodically forced pendulum (\ref{pendper}).
Such pseudo solutions turn $ k $ times along the separatrices
and are asymptotic to the equilibrium $0$, mod $ 2 \pi $, 
for $ t \to \pm \infty $.
More precisely in lemma 2.4 of \cite{BB3} it is proved that 
for all $ k \in {\bf N} $, for all 
$ \teta_1 < \ldots < \teta_k $ with
$ \min_i ( \teta_{i+1} - \teta_i ) > L $, with 
$ L $ sufficiently large, independently of $ \om_\e $ and $ \mu $,
there exists a unique pseudo-homoclinic solution
$ q_{A, \teta }^L ( t ): {\bf R} \to {\bf R} $
which is a true solution of (\ref{pendper})
in each interval $ (- \infty, \teta_1) $, $ ( \teta_i, \teta_{i+1} ) $
($ i = 1, \ldots, k-1 $), $ ( \teta_k, +\infty ) $ and
$ q_{A, \teta }^L ( \teta_i ) = \pi ( 2 i - 1 ) $,
$ q_{A, \teta }^L (t) = q^\mu_{ A, \teta_1 } ( t ) $ in
$( - \infty, \teta_1 ) $ and
$ q_{A, \teta }^L ( t ) = 2 \pi k + q^\mu_{A, \teta_k }(t) $ in
$( \teta_k, + \infty )$.
Such pseudo-homoclinic orbits are found via the Contraction Mapping Theorem,
as small perturbations of a chain of ``1-bump homoclinic solutions''
$ q^\mu_{A,\teta_i} $. 

Then we consider the Lagrangian action functional evaluated
on these pseudo-homoclinic  orbits  $ q_{ A, \teta }^L $ 
depending on $ n + k $ variables
$$
F_\mu^k (A_1, \ldots, A_n, \teta_1, \ldots, \teta_k ) =
\int_{- \infty}^{+\infty} {\cal L}_{\mu, A} (q_{A,\teta}^L (t),
{\dot{  q}}_{A,\teta}^L (t),t) \ dt.
$$
Setting $ e_k = (1, \ldots , 1 ) \in {\bf R}^k $, the following invariance
property, inherited from the autonomy of ${\cal H}_\mu$,  holds
\be \label{invariance}
F_\mu^k ( A, \theta + \eta e_k ) = F_\mu^k ( A + \eta \om_\e, \teta ),
\qquad  \forall \teta \in {\bf R}^k, \eta \in {\bf R}.
\ee
Let $ {\cal F}_\mu^k : {\bf T}^n \times {\bf R}^k \to {\bf R} $
be the  ``$ k $-bump heteroclinic function'' defined by
\be\label{eq:kbhf}
{\cal F}_\mu^k ( A , \teta ) := F_\mu^k (A, \teta) - (I_0'- I_0) \cdot A.
\ee
\begin{lemma}\label{lem:heter}
$\forall I_0, I_0' \in {\bf R}^n $,
if $ ( A, \teta ) $ is a critical point of the 
``$ k$-bump heteroclinic function''
$ {\cal F}_\mu^k (A, \teta) $, then
$ ( I_\mu (t), \om_\e t + A,$ $ q_{A, \teta}^L (t),
{\dot q}_{A, \teta }^L (t) )$
where
$ I_\mu (t) = I_0 - \mu \int_{-\infty}^t 
(\cos q_{ A, \teta }^L ( s ) - 1 ) \partial_\vphi f (\om_\e s+ A)ds $
is an heteroclinic solution connecting  $ {\cal T}_{I_0} $ to
$ {\cal T}_{I_0'} $.
\end{lemma}

By lemma \ref{lem:heter}, in order to get heteroclinic
solutions connecting ${\cal T}_{I_0}$ to ${\cal T}_{I_0'}$, 
we need to find
critical points of $ {\cal F}_\mu^k (A, \teta) $.
When $\min_i (\teta_{i+1} - \teta_i) \to + \infty $
the ``$k$-bump homoclinic function'' $ F_\mu^k ( A, \teta) $ turns out to be
well approximated simply by the sum of the functions
$ F_\mu ( A, \teta_i) $ according to the following lemma.
We set $ \teta_0 = - \infty$ and $ \teta_{k+1} = + \infty $.

\begin{lemma}\label{approxsum}
There exist positive constants $ C_1, L_1 > 0 $ and functions
$  R_i (\mu, A, \teta_{i-1}, \theta_i, \theta_{i+1}) $ such that
$ \forall \e >0 $, $ \forall  0<  \mu < \mu_0 $,
$ \forall L > L_1 $, $ \forall \teta_1 < \ldots < \teta_k $ with
$ \min_i ( \teta_{i+1} - \teta_i ) > L $
\be\label{eq:somma}
F^k_\mu (A, \theta_1, \cdots,
\theta_k ) = \sum_{i=1}^k F_\mu (A, \theta_i )+
\sum_{i=1}^k R_i (\mu, A, \teta_{i-1}, \theta_i, \theta_{i+1}),
\ee
with
\be\label{eq:resti}
|R_i (\mu, A, \teta_{i-1}, \theta_i, \theta_{i+1}) |
\leq C_1 \exp (- C_1 L ). 
\ee
\end{lemma}

\section{The splitting condition}

We now give a general ``splitting condition'' on the homoclinic
function $ G_\mu $
well suited to describe the non-uniform splitting of stable 
and unstable manifolds
which takes place in three time scale systems. 
Roughly, the ``splitting condition'' \ref{cond} below
states that $ G_\mu $ possesses a maximum and provides explicit 
estimates of the non-uniform splitting. It will be used, 
in the next section, to prove the 
shadowing theorem \ref{step1}.
As a paradigmatic example, 
we will verify, in lemma \ref{step3}, that, when the perturbation 
$ f( \vphi ) =  \sum_{j=1}^n \cos \vphi_j$, 
the ``splitting condition'' is satisfied,
see also remark \ref{rem:moregen}. 

\begin{condition}\label{cond} \ {\bf ``Splitting Condition''.}
There exist $ \ov{A} \in {\bf R}^n $ 
and a basis $\{ \Omega_1, \ldots, \Omega_n \} $ of $ {\bf R}^n $,
$ n \geq 3 $,
such that $ \om_\e \in {\bf R}_+ \Omega_1$, $ 1/2 \leq | \Omega_i| \leq 2 $,
${\rm det} \{ \Omega_1,\ldots,\Omega_n \} \geq 1/2 $,
$ \{ \Omega_3, \ldots, \Omega_n \}$ is 
an orthonormal basis of $\{ \Omega_1, \Omega_2 \}^{\bot}$,
and which enjoy the following properties : let us
define $ H_\mu (a_1 , \ldots, a_n ) $ as the homoclinic
function $ G_\mu (A) $ in the new basis, namely 
\be
H_\mu (a_1,\ldots, a_n) = G_\mu (\ov{A} + a_1\Omega_1+ \ldots + a_n \Omega_n).
\ee
Then there exist positive constants $ \rho, \s , \d_1, \d_2, \d_3  >0$,
with $ 3 \s < \rho $, $ \d_2 < \d_3 $, and two
continuous functions $ l_1, l_2: [- \rho, \rho ] 
\times \ov{B}_\rho^{n-2} \to {\bf R}$
with $ l_1 (x ) < l_2 (x ) $ for all 
$ x \in [ - \rho, \rho ] \times \ov{B}_{\rho}^{n-2}$, such that:
\begin{itemize}
\item $(i)$ for  $ x = (a_2, \ldots, a_n) 
\in [-\rho,\rho] \times\ov{B}_{\rho}^{n-2}$ 
$$
{\cal J} ( x ) := \sup_{a_1 \in [ l_1(x), l_2(x) ]}
H_\mu (a_1,x) \geq \max \Big\{ H_{\mu}(l_1(x),x),
H_{\mu}(l_2(x),x) \Big\} + \d_1;
$$
\item $(ii)$ for all $y=(a_3,\ldots ,a_n) \in \ov{B}_{\rho}^{n-2}$,
$$
\forall a_2 \in [-\s,\s], \  \  {\cal J}(a_2,y) \geq {\cal
J}(0,y) - \frac{\d_2}{2},
$$
$$
\forall a_2 \in [-\rho, -\rho+2\s] \cup [\rho-2\s,\rho], \  \
{\cal J}(a_2,y) \leq {\cal J}(0,y)-\d_2;
$$
\item $(iii)$ 
$$
\forall a_2 \in [-\s,\s], \ \ \forall y \in \ov{B}_{\s}^{n-2}, \ 
 \ {\cal J}(a_2,y) \geq {\cal J}(0,0) -\frac{\d_3}{2},
$$
$$
\forall a_2 \in [-\rho,\rho], \  \ 
\forall y \in  \ov{B}_{\rho}^{n-2} \backslash B_{\rho-2\s}^{n-2}, \  \
{\cal J}(a_2,y) \leq {\cal J}(0,0)-\d_3.
$$
\end{itemize}
\end{condition}

The next lemma states that the former ``splitting condition''
is satisfied by the homoclinic function $ G_\mu $ if (and only if) 
it holds for the homoclinic function $ \wtilde{G}_\mu $.

\begin{lemma} \label{step2}
Assume that $ \wtilde{G}_\mu $ satisfies the splitting condition \ref{cond}
with maps $ \wtilde{l}_1, \wtilde{l}_2 $ and parameters 
$ \rho $, $ \sigma $, $ \d_1 $, $ \d_2 $, $ \d_3 $. 
Then $ G_\mu $ satisfies the splitting condition \ref{cond} as well,
for some maps $ l_{1,2} = \wtilde{l}_{1,2} + O(\mu / \sqrt{\e})$ and 
with the same parameters. The converse is also true. 
\end{lemma}

\begin{pf}
By theorem \ref{thm:diffeo}, 
$ {\wtilde G}_\mu = G_\mu \circ \psi_\mu $,
where $ \psi_\mu ( A ) = A + k_\mu ( A ) \om_\e $ and $ \psi_\mu $ is a
homeomorphism. Set $ \wtilde{H}_\mu ( a_1, \ldots , a_n ) = 
\wtilde{G}_\mu (\ov{ A } + a_1 \Omega_1 + \ldots + a_n \Omega_n )$. 
We have
$$
\wtilde{H}_\mu ( a_1, a_2,  \ldots , a_n ) = 
H_\mu \Big( a_1 + \ov{k}_\mu
( a_1, \ldots, a_n ) \frac{ |\om_\e|}{| \Omega_1 | }, a_2, \ldots , a_n \Big),
$$
where $ \ov{k}_\mu ( a_1, \ldots, a_n) := k_\mu (\ov{A} + a_1 \Omega_1 +
\ldots + a_n \Omega_n)$.

Assume that $ \wtilde{G}_\mu $ satisfies condition \ref{cond}
with maps $ \wtilde{l}_1,\wtilde{l}_2$.  
For all $ x = ( a_2, \ldots, a_n) \in [-\rho,\rho] \times \ov{B}_{\rho}^{n-2}$, 
the map $ a_1 \mapsto a_1 + \ov{k}_{\mu} (a_1, x ) |\om_\e|/ |\Omega_1| $ is a
homeomorphism from the interval $ ( \wtilde{l}_1 (x),
\wtilde{l}_2(x))$ to the interval $(l_1(x),l_2(x))$, 
where $ l_j ( x ) : = \wtilde{l}_j (x) + \ov{k}_\mu
(\wtilde{l}_j(x), x) |\om_\e| / |\Omega_1| $ ($ j = 1,2 $). 
There results that, for all $ x = ( a_2, \ldots, a_n) \in$ $
[-\rho,\rho] \times \ov{B}_{\rho}^{n-2}$ 
$$ 
\wtilde{\cal J} (x) := \sup_{a_1 \in 
[ \wtilde{l}_1 ( x ), \wtilde{l}_2 ( x ) ]}
{\wtilde H}_\mu (a_1,x) = \sup_{a_1 \in 
[l_1 ( x ), l_2 ( x ) ]} H_\mu (a_1,x) = { \cal J} (x). 
$$ 
Therefore $ G_\mu $ satisfies the splitting condition
\ref{cond}, with maps $ \wtilde{l}_j $ replaced by $ l_j $, and
the same positive parameters.
Since $ \ov{k}_\mu = O(\mu) $ and $ | \om_\e | = O( 1 / \sqrt{\e} ) $
we have $ | l_j ( x ) - \wtilde{l}_j (x) | = O( \mu / \sqrt{\e}) $.
\end{pf}

We now give a paradigmatic example where the former ``splitting condition''
is satisfied. Assume that the perturbation $f$ is given by  
$ f(\vphi_1, \ldots, \vphi_n ) =  \sum_{j=1}^n \cos \vphi_j . $
In the next lemma we show that the corresponding 
homoclinic function $ \wtilde{G}_\mu $ 
satisfies the ``splitting condition'' \ref{cond} and hence, 
by lemma \ref{step2}, $ G_\mu $ as well satisfies 
the ``splitting condition'' \ref{cond}.

\begin{lemma} \label{step3} 
Assume that $ f( \vphi ) = \sum_{j=1}^n \cos \vphi_j $. There 
exist a basis $\{ \Omega_1, \ldots, \Omega_n \}$ and 
a positive constant $ \d_0 $ such that, if $ \e $ is small,
 $ 0 <  \mu \e^{-3/2} < \d_0 $ and 
$ 0 < \mu \e^{-2a-1} < \d_0 $, 
then $ \wtilde{G}_\mu $ satisfies 
the ``splitting condition'' \ref{cond}, with $\ov {A} = 0$,  
$ \rho= \pi \e^{a+1/2} $,
$ \sigma = \rho / 6 $, $ \d_1 = \d_3 = \mu \rho^2/2 $, 
$ \d_2 = 3\pi \mu \e^{-1/2} \exp (-\pi/(2\sqrt{\e}))$, 
$ \wtilde{l}_1 (x) = - 2 \pi, \wtilde{l}_2 (x) = 2 \pi $.
\end{lemma}

\begin{pf}
In order to simplify the notations we give the proof for $ n = 3 $ and
we assume that $ | \beta | = 1 $. 
We will prove that $ \wtilde{G}_\mu $ satisfies 
the ``splitting condition'' \ref{cond} with 
$\ov{A}=0$ and w.r.t the basis
$$
\Omega_1=(1,\e^{a+1/2} \beta), \quad  \quad 
\Omega_2=(0,\beta), \quad \quad 
\Omega_3=(0,\beta'),
$$
where $ | \beta' | = 1 $ and $ \beta \cdot \beta' =0 $.
We set $\rho =\pi \e^{a+1/2}$ and we assume that 
$ 0 < \mu \leq \delta \rho^2 $, $0 < \mu \leq \d \e^{3/2} $, 
where $ \d $ is a small constant (independent of $ \e $) to be specified later.
Let $\ov{\delta} >0$ be such that theorem  \ref{thm:tts}  holds
for $0 <\mu \leq \ov{\d} \e^{3/2}$. We shall always choose 
$0 < \d \leq \ov{\d}$. 

From now on, notation $K_i$ will be used for positive universal
constants, whereas notation $ c_i(\d) $ will be used for
positive constants depending only on $ \d $. Notation 
$ u = O(v) $ will mean that there exists a universal constant
$ K $ such that $|u| \leq K|v|$. 
\\[1mm]
\indent
Our first aim is to prove expression (\ref{Hmuap}) below.
It easily results that, if $ f( \vphi ) = \sum_{j=1}^3 \cos \vphi_j $,
\be\label{casopart}
\Gamma_0 (\e, A_2 ) = \sum_{j=2}^3
\frac{ 2 \pi \b_j \e^a }{ {\rm sinh} (\b_j \e^a \frac{\pi}{2})} \cos A_j  
\qquad {\rm and} \qquad 
\Gamma_1 (\e, A _2) = 
\frac{ \pi}{ \sqrt{\e} {\rm sinh} (\frac{\pi}{2 \sqrt{\e}})}.
\ee
By thereom \ref{thm:tts} we have
\be \label{Gmu1}
\wtilde{G}_{\mu}(A_1,A_2,A_3)=\wtilde{g}_0(A_2,A_3)+ 2 {\rm Re} \ 
\Big[ \wtilde{g}_1 (A_2,A_3) e^{iA_1}\Big] + 
O \Big( \mu \e^{-1/2} e^{-\pi/\sqrt{\e}} \Big)
\ee
and, by (\ref{casopart}), up to a constant that we shall omit, 
\be \label{G0}
\wtilde{g}_0(A_2,A_3)= \frac{\mu 2\pi \beta_2 \e^a}{\sinh(\beta_2 \e^a
\frac{\pi}{2})} \cos A_2 + \frac{\mu 2 \pi \beta_3
\e^a}{\sinh(\beta_3 \e^a
\frac{\pi}{2})} \cos A_3  + O(\mu^2),
\ee
\be \label{G1}
\wtilde{g}_1(A_2,A_3)= \frac{\mu \pi}{\sqrt{\e}\sinh(
\frac{\pi}{2\sqrt{\e}})} + O \Big( \frac{\mu^2}{\e^2}
e^{-\pi/2\sqrt{\e}} \Big).
\ee
In this proof we shall use the abbreviations
$$
C_\e=\frac{2 \pi \beta_2 \e^a}{\sinh(\beta_2 \e^a
\frac{\pi}{2})}+ \frac{2 \pi \beta_3 \e^a}{\sinh(\beta_3 \e^a
\frac{\pi}{2})}, \quad  \quad D_{\e}= \frac{2 \pi}{\sqrt{\e}\sinh(
\frac{\pi}{2\sqrt{\e}})}.
$$
Note that, as $ \e \to 0 $, we have 
\be \label{corrr}
\frac{2 \pi \beta_j \e^a}{\sinh(\beta_j \e^a
\frac{\pi}{2})}=4+O(\e^a), \quad  \quad D_\e=\frac{4\pi}{\sqrt{\e}}
e^{-\pi/(2\sqrt{\e})} \Big( 1+O(e^{-\pi/\sqrt{\e}}) \Big).
\ee
We shall consider $ \e $ small so that
\be \label{estCeDe}
 \frac{3\pi}{\sqrt{\e}}
 e^{-\pi/2\sqrt{\e}} \leq D_\e \leq \frac{5\pi}{\sqrt{\e}}
 e^{-\pi/2\sqrt{\e}}
\ee 
By (\ref{Gmu1}) and (\ref{G1}), 
since $0 < \mu \leq \d \e^{3/2}$,
\be \label{expGmu}
\wtilde{G}_{\mu}(A_1,A_2,A_3)=\wtilde{g}_0(A_2,A_3)+ \mu D_\e \cos A_1 +
O \Big( \mu \e^{-1/2} e^{-\pi/\sqrt{\e}} + 
\mu \d \e^{-1/2} e^{-\pi/2\sqrt{\e}} \Big).
\ee
Since $ ( A_1, A_2, A_3 ) = a_1 \Omega_1 + a_2 \Omega_2 + a_3\Omega_3 = 
(a_1, (a_1 \e^{a+1/2} +a_2) \beta + a_3 \beta' )$, 
the homoclinic function $ \wtilde{G}_\mu $ in the new basis 
$\{ \Omega_1, \Omega_2, \Omega_3 \} $ writes
\be\label{eq:primaexp}
\wtilde{H}_{\mu}(a_1,a_2,a_3)=\wtilde{G}_\mu 
\Big( a_1\Omega_1 + a_2 \Omega_2 + a_3
\Omega_3 \Big) = \wtilde{G}_\mu \Big( a_1, (a_1 \e^{a+1/2} +a_2) \beta +
a_3 \beta' \Big).
\ee
Define $ \wtilde{h}_0 (b_2,a_3) = \wtilde{g}_0 (b_2 \beta + a_3 \beta')$.
By (\ref{eq:primaexp}) and (\ref{expGmu}), there 
exists $ c_0 ( \d  ) > 0 $ such that, 
for all $ 0 < \e \leq c_0 ( \d ) $, 
\be \label{Hmu1}
\wtilde{H}_{\mu}(a_1,a_2,a_3)=
\wtilde{h}_0 \Big( a_1 \e^{a+1/2} +a_2,a_3 \Big) + \mu D_\e \cos a_1 +
O \Big( \mu \d \e^{-1/2} e^{-\pi/2\sqrt{\e}} \Big).
\ee
We derive from this latter expression and (\ref{estCeDe}) that  
\be \label{Hmu2}
\wtilde{H}_{\mu}(a_1,a_2,a_3)=
\wtilde{h}_0 \Big( a_1 \e^{a+1/2} +a_2,a_3 \Big) +
O \Big( \mu \e^{-1/2} e^{-\pi/2\sqrt{\e}} \Big).
\ee
By (\ref{G0}) and (\ref{corrr})
\be \label{G0ap}
\wtilde{g}_0(A_2,A_3)=\mu C_\e - 2\mu ( A_2^2 + A_3^2) 
+ O\Big(\mu \e^a (A_2^2 + A_3^2)\Big)+ 
O \Big( \mu (A_2^4 +A_3^4) \Big) + O(\mu^2)  .
\ee
We shall assume 
in the sequel of the proof that $ a_2, a_3 \in [-\rho,\rho]$,
$a_1 \in [-2\pi,2\pi]$, so that, since $ \rho = \pi \e^{a + 1/2} $, there
results $ a_1 \e^{a+1/2} \in [-2\rho , 2\rho]$, $ b_2 = a_1 \e^{a+1/2}+
a_2 \in [-3\rho,3\rho]$ and $b_2^4+a_3^4=O(\rho^4)$.
Moreover we have that $\mu^2 \leq \mu \d \rho^2$ and there
exists $ c_1(\d) \in (0,c_0(\d))$ such that, 
if $0 < \e \leq c_1(\d)$, then $\e^a \leq \delta$ and $\rho^4 \leq \d \rho^2$.   
Note also that, since $\b, \b'$ are orthonormal vectors, we have 
$ A_2^2 + A_3^2 = b_2^2 + a_3^2 $. Finally we derive 
from (\ref{G0ap}) that, for $ 0 < \e \leq c_1 ( \d )$, 
\be \label{H0}
\wtilde{h}_0 ( b_2, a_3) = \mu C_\e -2\mu ( 
b_2^2+a_3^2) +  O( \mu \d \rho^2 ).
\ee
Since $ \rho = \pi \e^{ a + 1/2 }$ we have  
$ \e^{-1/2} e^{-\pi / 2 \sqrt{\e}} = o ( \rho^2 )$ as 
$ \e \to 0$; therefore, by (\ref{Hmu2}) and (\ref{H0}), there exist 
$ K_0 > 0 $, $ c_2 (\d) \in (0, c_1(\d))$ such that,
for all $ 0< \e \leq c_2(\d)$, 
\be \label{Hmuap}
\wtilde{H}_\mu (a_1,a_2,a_3)=\mu C_\e - 2\mu
(b_2^2+a_3^2) + r_0(a_1,a_2,a_3), \quad
| r_0(a_1,a_2,a_3)| \leq K_0 \mu \d \rho^2 ,
\ee
where $ b_2 = a_1 \e^{a+1/2} + a_2$.
\\[1mm]
\indent
We now prove that point $(i)$ of the ``splitting condition'' \ref{cond}
is satisfied by $\wtilde{G}_{\mu}$ with $ \d_1=\mu \rho^2/2 $, 
$ \wtilde{l}_1(x) = -2\pi $ and 
$ \wtilde{l}_2 ( x ) = 2 \pi $ where $x := (a_2, a_3)$.
Let us consider 
$ {\cal J}(a_2,a_3) := \sup_{a_1\in [-2\pi,2\pi]} 
\wtilde{H}_{\mu} (a_1,a_2,a_3)$.
Since $ a_2 \in [ - \rho, \rho ] $, 
$ -a_2 \e^{-(a+1/2)} \in [- \pi, \pi ]$ and we can derive from (\ref{Hmuap})
that  
\be \label{minorJ}
{\cal J} (a_2,a_3) \geq \wtilde{H}_\mu \Big( -a_2 \e^{-(a+1/2)},a_2,a_3 \Big)
\geq \mu C_\e-2\mu  a_3^2 - K_0 \mu \d \rho^2 .
\ee
If $ a_1 = \pm 2 \pi $ then $ b_2 = a_1 \e^{a+1/2}
+ a_2 = a_2 \pm 2 \rho $ and then, since $a_2 \in [- \rho, \rho]$, we get
$ |b_2| \geq \rho $.
As a consequence, by (\ref{Hmuap}) and (\ref{minorJ}),
\be \label{Hmuestsup}
\wtilde{H}_\mu \Big( \pm 2\pi,a_2,a_3 \Big) \leq  \mu C_{\e} 
-2\mu ( \rho^2+a_3^2) + K_0 \mu \d \rho^2  \leq
{\cal J}(a_2,a_3)-2\mu \rho^2 + 2K_0 \mu \d \rho^2 .
\ee
Choosing  $ \d < 1 / 2 K_0 $, we get in (\ref{Hmuestsup}) that 
$\wtilde{H}_\mu (\pm 2\pi,a_2,a_3)\leq 
{\cal J}(a_2,a_3)-\mu \rho^2/2. $
It results that condition \ref{cond}-$(i)$ 
is satisfied with $ \d_1 = \mu \rho^2/2 $,
$ l_1 ( x ) = - 2 \pi $ and $ l_2 (x) = 2 \pi $ 
where $ x = (a_2, a_3) $.
\\[2mm]
\indent  We now turn to the proof of $(ii)$ and $(iii)$.
If $a_2,a_3 \in [-\rho,\rho]$, $a_1 \in [-2\pi,2\pi]$ and  
$| b_2 | = |a_1 \e^{a+1/2}+a_2| \geq \sqrt{2K_0\d}\rho$,  
then, by (\ref{Hmuap}) and (\ref{minorJ}),  
$$
\wtilde{H}_\mu (a_1,a_2,a_3) \leq \mu C_\e - 2\mu
( a_3^2+ 2K_0 \d \rho^2 ) + 
K_0 \mu \d  \rho^2  \leq \mu C_\e - 2\mu
 a_3^2 - 3K_0 \mu \d \rho^2 < {\cal J}(a_2,a_3).
$$
Hence
\be \label{expJ}
{\cal J} ( a_2, a_3) = \sup \Big\{ \wtilde{H}_\mu (a_1,a_2,a_3) \ ;
\ a_1 \e^{a+1/2} \in 
\Big[ -a_2-\sqrt{2K_0\d}\rho, -a_2+\sqrt{2K_0\d} \rho \Big] \Big\}.
\ee
We use here that, since $ 2 K_0 \d < 1 $, 
$ [ - a_2 - \sqrt{2K_0\d} \rho, -a_2 + \sqrt{2K_0\d} \rho ]
\subset [-2\rho, 2\rho]=[-2\pi \e^{a+1/2} , 2\pi \e^{a+1/2}]$. 
Writing $ a_1 = (b_2 - a_2) \e^{- (a + 1/2)} $, we derive from 
(\ref{expJ}) and (\ref{Hmu1}) that 
\begin{eqnarray} \label{eq:3.17}
{\cal J}(a_2,a_3) & = & \sup_{b_2\in 
[-\sqrt{2K_0\d}\rho,\sqrt{2K_0\d}\rho]} \wtilde{H}_{\mu}
 \Big( (b_2-a_2)\e^{-a-1/2},a_2,a_3 \Big) \nonumber \\ 
& = & \sup_{b_2 \in [-\sqrt{2K_0\d}\rho,\sqrt{2K_0\d}\rho] }
\Big(\wtilde{h}_0 (b_2,a_3) +\mu D_\e \cos 
\Big( \dps \frac{b_2-a_2}{\e^{a+1/2}} \Big)
\Big) + O \Big( \delta \mu \e^{-1/2} e^{-\pi/2\sqrt{\e}} \Big).
\end{eqnarray}
Now, if $b_2 \in [-\sqrt{2K_0\d}\rho,\sqrt{2K_0\d}\rho]$ then
$b_2 \e^{-a-1/2} \in [-\pi\sqrt{2K_0\d},\pi\sqrt{2K_0\d}]$, so we
can write that
\be\label{cosexp}
\cos \Big( \frac{b_2-a_2}{\e^{a+1/2}} \Big) 
= \cos \Big( \frac{-a_2}{\e^{a+1/2}} \Big) + O(\sqrt{\d}).
\ee
As a consequence, by (\ref{eq:3.17}) and (\ref{estCeDe}) there holds
\begin{eqnarray}\label{exprJ}
{\cal J}(a_2,a_3) 
& = & \sup_{b_2 \in [-\sqrt{2K_0\d}\rho,\sqrt{2K_0\d}\rho] }
\Big(\wtilde{h}_0 (b_2,a_3) +
\mu D_\e \cos \Big( \frac{-a_2}{\e^{a+1/2}} \Big)
\Big) + O \Big( \sqrt{\delta} \mu \e^{-1/2} e^{-\pi/2\sqrt{\e}} \Big)
\nonumber \\ & = & 
\wtilde{m}(a_3) +\mu D_\e \cos \Big( \frac{a_2}{\e^{a+1/2}} \Big)
+ O \Big( \sqrt{\delta} \mu \e^{-1/2} e^{-\pi/2\sqrt{\e}} \Big),
\end{eqnarray}
\noindent
where we have set
\be\label{def:mtilde}
\wtilde{m}(a_3):= \sup_{b_2 \in [-\sqrt{2K_0\d}\rho,\sqrt{2K_0\d}\rho]}
\wtilde{h}_0(b_2,a_3).
\ee
Finally, there exists $ K_1>0 $ such that,
by (\ref{exprJ}), 
\be \label{demii}
{\cal J}(a_2,a_3)={\cal J}(0,a_3)+ 
\mu D_\e \Big( \cos \Big( \frac{a_2}{\e^{a+1/2}} \Big) - 1 \Big) +
r_1(a_2,a_3), \quad  |r_1(a_2,a_3)| \leq K_1 
 \sqrt{\delta}\mu \e^{-1/2} e^{-\pi/2\sqrt{\e}}.
\ee
We are now in position to prove condition \ref{cond}-$(ii)$. Assume 
$ 0 < \d \leq \pi^2 / 4K_1^2$ and choose $ \s = \rho/6 =
\e^{a+1/2} \pi /6$ and  $\d_2= 3\pi \mu \e^{-1/2} e^{-\pi/2\sqrt{\e}}$.
If $ a_2 \in [- \s, \s] $ then 
$ \cos (a_2/ \e^{a+1/2})-1 \geq -1+\sqrt{3}/2 \geq -1/6$. This
readily implies, by (\ref{demii}) and (\ref{estCeDe}), that
$ { \cal J }( a_2, a_3 ) \geq {\cal J}(0, a_3)- \delta_2/2$.
If $a_2 \in [-\rho,-\rho+2\s]\cup
[\rho-2\s,\rho]$ then $a_2/\e^{a+1/2} \in [-\pi,-2\pi/3]\cup 
[2\pi/3,\pi]$, so that  
$ \cos (a_2/ \e^{a+1/2})-1  \leq -3/2$.  It follows,
still by (\ref{demii}) and (\ref{estCeDe}), that 
$ { \cal J}(a_2,a_3) \leq {\cal J}(0,a_3)- \delta_2 $.
This proves condition \ref{cond}-$(ii)$.
 
In order to prove condition \ref{cond}-$(iii)$, 
we notice that, by (\ref{H0}) and 
the definition of $\wtilde{m}$ given in (\ref{def:mtilde}),  
\be\label{expmtilde}
\wtilde{m}(a_3)=\mu C_{\e}-2\mu a_3^2 +
O(\delta \mu \rho^2). 
\ee
Hence there exist $ K_2 > 0 $ and $ c_3 ( \d ) \in (0,c_2(\d))$ such
that, for all $ 0 < \e \leq c_3 ( \d )$, by (\ref{expmtilde}), (\ref{exprJ})
and (\ref{estCeDe}) 
\be \label{demiii}
{\cal J}(a_2,a_3)= {\cal J}(0,0)- 2 \mu a_3^2  +
r_2(a_2,a_3), \quad {\rm with} \quad |r_2(a_2,a_3)|\leq K_2 \delta \mu \rho^2.
\ee
Let us assume $\delta \leq 1/ 6K_2 $ and let
$ \d_3 = \mu \rho^2/2 $.
By (\ref{demiii}) and (\ref{estCeDe}), if 
$a_3 \in \ov{B}_{\s}^{n-2}$, then ${\cal J}(a_2,a_3) \geq 
{\cal J}(0,0) -\d_3/2$; if $a_3 \in \ov{B}_{\rho}^{n-2}
\backslash {B}_{\rho-2\s}^{n-2}$ then 
${\cal J}(a_2,a_3) \leq {\cal J}(0,0) - \delta_3$.
\\[1mm]
\indent
As a conclusion, lemma \ref{step3} holds with
$\d_0=\min \{ \ov{\d}, 1/ 2K_0 , \pi^2/4K_1^2, 1/ 6 K_2 \}$.
\end{pf}

\begin{remark}\label{rem:moregen}
The former splitting condition holds also also for more general 
perturbations $ f( \vphi_1, \ldots, \vphi_n ) $ 
for which $ f_0 ( \vphi_2, \ldots, \vphi_n ) $ possesses a nondegenerate
maximum at $ (\ov{\vphi}_2, \ldots, \ov{\vphi}_n )$ 
and $ f_1 (\ov{\vphi}_2, \ldots, \ov{\vphi}_n ) \neq 0 $ where  
$ f_{k_1} ( \vphi_2, \ldots, \vphi_n) = (1 / 2 \pi ) 
\int_0^{2 \pi}
f( \sigma, \vphi_2, \ldots, \vphi_n ) e^{-i k_1 \sigma} \ d \sigma.$
This kind of condition is considered in theorem 5.2 of \cite{BB3}.
\end{remark}

\section{The shadowing theorem}

In this section we shall prove, under the ``splitting condition''
\ref{cond}, our general shadowing theorem.

\begin{theorem} \label{step1}
Let $ n \geq 3 $ and assume that the homoclinic function 
$ G_\mu $ satisfies the splitting condition
\ref{cond}. Let $ \om_\e $ be a $(\gamma_\e, \tau)$-diophantine vector, i.e.
$ | \om_\e \cdot k | \geq \gamma_\e / |k|^\tau $
$\forall k \in {\bf Z}^n \backslash \{ 0 \}$.
Then, for all $I_0, I'_0 \in {\bf R}^n $ such that 
$( I_0' - I_0 ) \in $
${\rm Span} \{ \Omega_3,\ldots,\Omega_n \}$, 
there exists an heteroclinic trajectory from 
$ {\cal T}_{I_0}$ to ${\cal T}_{I_0'}$
which connects a $ \eta$-neighbourhood of torus ${\cal T}_{I_0}$ to 
a $ \eta$-neighbourhood  of torus ${\cal T}_{I'_0}$ in the 
``diffusion time''
\be\label{timediff}
T_d \leq 
C \frac{\rho |I'_0-I_0|}{\d_3} \ {\rm max} 
\Big\{ \frac{1}{ \g_\e \s^\tau}, |\ln \d_1|, 
|\ln \d_2|, \frac{\Delta}{| \om_\e |} \Big\}   
+ |\ln \eta|, 
\ee
where 
$ \Delta := \Big\{ 
\max_{x \in [-\rho,\rho] \times\ov{B}_{\rho}^{n-2} } l_2 (x) - 
\min_{x \in [-\rho,\rho] \times\ov{B}_{\rho}^{n-2} } l_1 (x ) \Big\} $.
\end{theorem}

\begin{remark}
The diophantine condition on the frequency vector $ \om_\e $ 
restricts the values of $ \e $ and $\beta$ that we consider.
In any case, if for instance $\beta$ is
($\gamma$,$n-2$)-diophantine then   for $ \tau \geq n - 1 $
there exist $c_0 >0$ and a sequence $ \e_j \to 0 $ such that 
$\om_\e $ is $( \gamma_\e,\t )$-diophantine
with  $ \g_\e = c_0 \e^a $, see for example \cite{GGM}.
\end{remark}

\begin{remark}
The meaning of (\ref{timediff}) is the following:
the diffusion time $ T_d $ is estimated by the product
of the number of heteroclinic transitions
$ k = $ ( heteroclinic jump / splitting ) $ = | I_0' - I_0 | / \delta_3 $,
and of  the time $ T_s $ required for a single transition,
that is $ T_d \approx k \cdot T_s  $.
The time for a single transition $ T_s $ is bounded by the maximum
time between the
``ergodization time'' $ ( 1 / \gamma_\e \s^\tau ) $, i.e.
the time needed for the flow $ \om t $ to make an
$ \s $-net of the torus,
and the time $ \max \{ | \ln \delta_1 |,| \ln \delta_2 |, 
\Delta / | \om_\e | \}$ 
needed  to ``shadow''
homoclinic orbits for the forced pendulum equation. We use
here that these homoclinic orbits are exponentially
asymptotic to the equilibrium.

We could prove also the existence of connecting 
orbits for all 
$I_0' - I_0 \in $ Span $\{ \Omega_2,\ldots,\Omega_n \}$. 
In this case the number $ k $ of heteroclinic transitions would depend
also on $ \d_2 $,  see remark \ref{poliexp}. 
\end{remark}

\begin{pf}
Still for simplicity of notation we write the proof for $ n = 3 $. 
Then $ I'_0-I_0 = \pm |I'_0 -I_0| \Omega_3 $; 
we assume for definitiveness that $ I'_0-I_0 = |I'_0 -I_0| \Omega_3$, so that
$ (I'_0-I_0) \cdot (\sum_{j=1}^3 a_j \Omega_j)=|I'_0-I_0| a_3$.

We choose the number of heteroclinic transitions as
\be\label{numberhet}
k = \Big[ \frac{ 8 | I_0'-I_0 | \rho }{ \d_3 } \Big] + 1.
\ee
By lemma \ref{lem:heter}, in order to prove the theorem, it is sufficient
to find a critical point 
of the $ k $-bump heteroclinic function 
${\cal F }_\mu^k: {\bf T}^3 \times {\bf R}^k \to {\bf R}$
such that 
\be\label{eq:qtimed}
\theta_k - \theta_1 = O \Big( \frac{\rho |I'_0-I_0|}{\d_3} 
{\rm max} \Big\{ \frac{1}{\gamma_{\e} \s^{\tau}}, 
|\ln \d_1| , |\ln \d_2|, \frac{\Delta}{| \om_\e |}\Big\} \Big).
\ee 
We introduce suitable coordinates
$ (a_1, a_2 , a_3, s_1, \ldots, s_k ) \in {\bf R}^3 \times
(- 2 \pi, 2 \pi)^k $ defined by
\be\label{eq:cov}
A = \ov{A} + \sum_{j = 1}^3 a_j \Omega_j \quad {\rm and}
\quad
\forall i = 1, \ldots, k, \quad 
\teta_i = \frac{(\eta_i + s_i - a_1)|\Omega_1|}{|\om_\e |},   
\ee
where  $ \eta_i $ are constants to be chosen later.
Let $ H_{\mu}^k (a,s)= F_{\mu}^k (A,\theta) $
be the ``$k$-bump homoclinic function'' and 
$ {\cal H}_{\mu}^k (a,s) = {\cal F}_{\mu}^k (A,\theta)$
be the ``$k$-bump heteroclinic function'' expressed 
in the new variables $(a,s)$.

The function $ {\cal H}_\mu^k $ does not depend on $ a_1 $, since, 
by the invariance property (\ref{invariance}) (we 
recall that $ \Omega_1 = |\Omega_1| \omega_\e / | \omega_\e |$),
\begin{eqnarray}
{\cal H}_\mu^k (a,s) & = &
F_\mu^k \Big( \ov{A} + \sum_{j = 1}^3 a_j \Omega_j,
\frac{ (\eta_1 + s_1 - a_1)|\Omega_1|}{| \omega_\e |} , \ldots,
\frac{ (\eta_k + s_k - a_1)|\Omega_1|}{| \omega_\e |} \Big) 
- (I_0'-I_0) \cdot \Big( \sum_{j = 1}^3 a_j \Omega_j \Big) 
\nonumber \\
& = & F_\mu^k \Big( \ov{A} + \sum_{j = 2}^3 a_j \Omega_j ,
\frac{(\eta_1 + s_1)|\Omega_1|}{| \omega_\e |} , \ldots,
\frac{(\eta_k + s_k)|\Omega_1|}{| \omega_\e |} \Big) - 
|I_0'-I_0| a_3.
\end{eqnarray}
\noindent
In the sequel of the proof we shall use the abbreviation $ {\cal H}_{\mu}^k =
{\cal H}_\mu^k (0,a_2,a_3,s)$.
\\[1mm]
\indent 
We now choose the constants $(\eta_1, \ldots, \eta_k ) \in {\bf R}^k $.
Note that, since $ \om_\e $ is  $( \gamma_\e, \tau)$-diophantine, 
$ \Omega_1 $ satisfies the diophantine condition 
$$
|\Omega_1 \cdot k | \geq \frac{\g_\e |\Omega_1| }{ |\om_\e| |k|^\t}, \ \ 
\forall k \in {\bf Z}^n \backslash \{ 0 \}.
$$
Hence, by the results of \cite{BGW}, there exists $\ov{C} >0 $  
such that the ``ergodization time'' $ T_e $ of the torus $ {\bf T}^3 $
run by the linear flow $ \Omega_1 t $,  
i.e the smallest time for which $\{ \Omega_1 t \ ; \ 0\leq
t \leq T_e \}$
is a $\s-$ net of the torus, can be bounded from above by 
${\ov{C} |\om_\e |}/{ (\g_\e \s^\t)}$.
Hence for each interval $ J $ of length greater or equal to 
$ \ov{C} | \om_\e | / (\g_\e \s^\t) $ there exists $ \eta \in J$ 
such that 
\be \label{number}
d( \eta \Omega_1, 2 \pi {\bf Z}^3) < \s.
\ee
In particular there exists a constant $C_2$ and
there exist $ \eta_i $ such that
\be\label{eq:estisepa} 
\frac{ |\om_\e |}{C_1|\Omega_1|} \ln \Big(8 C_1 
\frac{1}{ \min \{ \d_1, \d_2 \} } \Big) + \Delta \leq
\eta_{i+1} - \eta_i \leq  \frac{ |\om_\e |}{C_1|\Omega_1|} 
\ln \Big(8C_1 \frac{1}{ \min \{ \d_1, \d_2 \}} \Big) +
\frac{C_2 |\om_\e |}{ \g_\e \s^\t} + \Delta, 
\ee
\be\label{theta2}
\eta_i \Omega_1 \equiv \chi_i, \ {\rm mod}  2 \pi {\bf Z}^3,
\quad \chi_i = y_i \Omega_2 + z_i \Omega_3 \quad {\rm with} \
|y_i| < \s, \ |z_i|<\s. 
\ee
\\[1mm]
\indent
In order to prove the theorem we just need to prove the existence
of a critical point of ${\cal H}_{\mu}^k$ in ${\bf R}^2
\times (\min l_1, \max l_2)^k $. 
The upperbound of the diffusion time 
given in (\ref{timediff}) will then be a consequence of
(\ref{eq:estisepa}) and (\ref{numberhet}). Indeed, 
by (\ref{eq:cov}) and (\ref{eq:estisepa}) we get that 
\be\label{eq:separ}
\teta_{i+1} - \teta_i = \frac{(\eta_{i+1} - \eta_i) |\Omega_1|}{|\om_\e |} 
+ \frac{( s_{i+1} - s_i) |\Omega_1|}{|\om_\e |} \leq 
\frac{1}{C_1} \ln \Big(8C_1 \frac{1}{ \min \{ \d_1, \d_2 \}} \Big) +
\frac{C_2 |\Omega_1 |}{ \g_\e \s^\t} + 
 \frac{ 2 \Delta | \Omega_1 | }{| \om_\e |}.
\ee
By (\ref{eq:separ}) there exists $ C >0 $ 
such that the time $\teta_{i+1} - \teta_i$ 
``spent for a single transition'' is bounded by 
\be\label{eq:singletra}
T_s := \max_i ( \teta_{i+1} - \teta_i ) \leq C 
{\rm max} \Big\{ \frac{1}{ \gamma_\e \s^\tau},|\ln \d_1| , |\ln \d_2|, 
\frac{\Delta}{|\om_\e|} \Big \}.
\ee
From (\ref{eq:singletra})
and (\ref{numberhet}) we derive immediately 
(\ref{eq:qtimed}) and then (\ref{timediff}).
\\[1mm]
\indent 
We now provide, using lemma \ref{approxsum}, a suitable expression of the
$k$-bump heteroclinic function ${\cal H}_\mu^k$. 
By lemma \ref{approxsum}, 
the invariance property (\ref{invpro}),
(\ref{theta2}) and since $ G_\mu : {\bf T}^3 \to {\bf R}$, we
get
\begin{eqnarray} \label{eq:espressionfi}
{\cal H}_\mu^k (a_2, a_3, s) & = &
\sum_{i = 1}^k \Big[ F_\mu \Big(\ov{A}+ \sum_{j = 2}^3 a_j \Omega_j,
\frac{ (\eta_i + s_i)|\Omega_1| }{ | \omega_\e |}  \Big) +
S_i (a_2,a_3,s_{i-1},s_i, s_{i+1}) - 
\frac{|I_0'-I_0|}{k} a_3 \Big] \nonumber \\
& = &
\sum_{i = 1}^k \Big[ G_\mu \Big(\ov{A}+ \sum_{j = 2}^3 a_j \Omega_j
+ \eta_i \Omega_1 + s_i \Omega_1  \Big) + 
S_i (a_2,a_3,s_{i-1},s_i, s_{i+1}) - 
\frac{|I_0'-I_0|}{k} a_3 \Big] \nonumber \\
& = & \sum_{i = 1}^k \Big[ H_\mu \Big( s_i , a_2+y_i , a_3+z_i
\Big) -\frac{|I_0'-I_0|}{k} a_3  + S_i \Big],   
\end{eqnarray}
where $ S_i := S_i (a_2,a_3,s_{i-1},s_i, s_{i+1}) =
R_i(A,\theta_{i-1},\theta_i,\theta_{i+1})$ 
after the change of variables (\ref{eq:cov}).
The left hand side inequality in (\ref{eq:estisepa}) implies that
$$
\theta_{i+1}-\theta_i \geq \frac{(\eta_{i+1}-\eta_i - \Delta) 
|\Omega_1|}{|\omega_{\e}|} \geq \frac{1}{C_1} \ln \Big(
8 C_1 \frac{1}{ \min \{ \d_1, \d_2 \} } \Big);
$$
 hence, by (\ref{eq:resti}),  
\be\label{eq:inteop}
| S_i | \leq \frac{\min \{ \d_1,\d_2 \}}{8}. 
\ee
\\[1mm]
\indent
We will maximize $ {\cal H}_\mu^k $ in the open set  
$$
U = \Big\{ (a_2,a_3,s) \in  {\bf R}^{k+2} \ \Big| \ 
\forall i \ \ a_2+y_i \in (-\rho,\rho), a_3+z_i \in (-\rho,\rho) , \ s_i \in 
\Big( l_1(a_2+y_i,a_3+z_i),
l_2(a_2+y_i, a_3+z_i) \Big) \Big\}.
$$
$ U \neq \emptyset $ since $ \{ 0 \} \times \{ 0 \}\times
\Pi_{i=1}^k (l_1(y_i,z_i), l_2(y_i, z_i)) \subset U$.
Since $U$ is bounded, ${\cal H}_{\mu}^k$ attains its maximum  over
$\ov{U}$ at some point $(\ov{a},\ov{s})=(\ov{a}_2,\ov{a}_3,\ov{s})$.
It is enough to prove that $(\ov{a},\ov{s}) \in U$. 
\begin{itemize}
\item We first prove that for all $i$, $\ov{s}_i \in
 (l_1(\ov{a}_2+y_i,\ov{a}_3+z_i), l_2(\ov{a}_2+y_i,
\ov{a}_3+z_i))$. Since  $(\ov{a},\ov{s})$ is a maximum point of
$ {\cal H}_\mu^k $ in $ \ov{U} $, 
for any $ t \in [ l_1(\ov{a}_2+y_i,\ov{a}_3+z_i), 
l_2 (\ov{a}_2+y_i,\ov{a}_3+z_i) ]$, replacing $\ov{s}_i$ with $t$
does not increase ${\cal H}_\mu^k$. Since such a substitution 
alters at most three terms among $ S_1, \ldots, S_k $
in (\ref{eq:espressionfi}), we obtain, 
using  (\ref{eq:inteop}), that for any $i$, for any $t\in
[l_1(\ov{a}_2+y_i,\ov{a}_3+z_i), 
l_2(\ov{a}_2+y_i,\ov{a}_3+z_i)]$,
$$
H_{\mu}(\ov{s}_i , \ov{a}_2 + y_i, \ov{a}_3 + z_i) \geq
H_{\mu}(t, \ov{a}_2 + y_i, \ov{a}_3 + z_i) - \frac{3}{4}\min \{ \d_1,\d_2 \}.
$$
Hence
$$
H_{\mu}(\ov{s}_i , \ov{a}_2 + y_i, \ov{a}_3 + z_i) \geq
{\cal J}(\ov{a}_2+y_i,\ov{a}_3+z_i) - \frac{3\d_1}{4},
$$
and, by condition \ref{cond}-$(i)$, this implies that
$\ov{s}_i \in (l_1(\ov{a}+\ov{\chi}_i),
l_2(\ov{a}+\ov{\chi}_i))$, where we have set $\ov{\chi}_i=(y_i,z_i)$.
\item We now prove that for all $i$, $\ov{a}_2+y_i \in (-\rho, \rho)$.
Indeed we have by (\ref{eq:espressionfi}) and  (\ref{eq:inteop})
\begin{eqnarray} \label{Hest}
{\cal H}_{\mu}^k ( \ov{a}_2 , \ov{a}_3 , \ov{s})& \leq &
\sum_{i=1}^k \Big[H_{\mu}(\ov{s}_i , \ov{a}_2 + y_i, \ov{a}_3 + z_i)
+ \frac{\min \{ \d_1,\d_2 \}}{8}-\frac{|I'_0-I_0|}{k} \ov{a}_3 \Big] 
\nonumber \\  
& \leq &
\sum_{i=1}^k \Big[{\cal J}( \ov{a}_2 + y_i, \ov{a}_3 + z_i)
+ \frac{\d_2}{8} -\frac{|I'_0-I_0|}{k} \ov{a}_3 \Big].
\end{eqnarray}
On the other hand, still by (\ref{eq:espressionfi}) and (\ref{eq:inteop}), 
choosing $ s = \wtilde{s} = ( \wtilde{s}_1, \ldots, \wtilde{s}_k ) $ with 
$ \wtilde {s}_i \in (l_1(y_i,\ov{a}_3+z_i), l_2 (y_i, \ov{a}_3+z_i))$
so that 
$ H_{\mu}( \wtilde{s}_i, y_i, \ov{a}_3+z_i)= {\cal J}(y_i,\ov{a}_3+z_i)$, 
we get
\begin{eqnarray}\label{4.8}
{\cal H}_\mu^k ( 0, \ov{a}_3 , \wtilde{s}) & = & 
\sum_{i = 1}^k \Big[ H_\mu \Big( \wtilde{s}_i , y_i , \ov{a}_3 + z_i\Big) 
-\frac{|I_0'-I_0|}{k} \ov{a}_3  + S_i \Big] \nonumber \\
& \geq & \sum_{i=1}^k \Big[{\cal J}(  y_i, \ov{a}_3 + z_i)
- \frac{\d_2}{8}- \frac{|I'_0-I_0|}{k} \ov{a}_3 \Big] \nonumber \\
& \geq &
\sum_{i=1}^k \Big[{\cal J}( 0, \ov{a}_3 + z_i)
- \frac{5\d_2}{8}- \frac{|I'_0-I_0|}{k} \ov{a}_3 \Big],
\end{eqnarray}
since $| y_i | < \s $ and by condition \ref{cond}-$(ii)$. Since 
$ {\cal H}_\mu^k (\ov{a},\ov{s}) \geq {\cal H}_\mu^k (0,\ov{a}_3,
\wtilde{s})$, we can  derive from  (\ref{Hest}) and (\ref{4.8}) that 
there exists $ i_0 \in \{ 1, \ldots, k \}$ 
such that $ {\cal J}( \ov{a}_2 + y_{i_0}, \ov{a}_3 + z_{i_0})
\geq {\cal J}( 0, \ov{a}_3 + z_{i_0})- 3\d_2/4$.  Still
by condition \ref{cond}-$(ii)$, 
$ \ov{a}_2 + y_{i_0} \in (- \rho + 2\s, \rho - 2 \s) $. As
a result, since (by (\ref{theta2})) 
 $|y_i - y_{i_0}| \leq 2 \s $ for all $ i $, we get that 
$\ov{a}_2 + y_i \in (-\rho,\rho)$. 
\item
At last we prove that, for all $ i $, $ \ov{a}_3 + z_i \in (-\rho,\rho) $.  
By (\ref{eq:espressionfi}) and (\ref{eq:inteop}), choosing 
$ s = \wtilde{s} = ( \wtilde{s}_1, \ldots, \wtilde{s}_k )$ such that 
$ \wtilde {s}_i \in (l_1(y_i,z_i), l_2 (y_i, z_i))$
so that  $H_\mu ( \wtilde{s}_i, y_i, z_i)= {\cal J}(y_i, z_i)$,
\be\label{4.9}
{\cal H}_{\mu}^k ( 0, 0 , \wtilde{s}) = 
\sum_{i = 1}^k \Big[ H_\mu \Big( \wtilde{s}_i , y_i ,  z_i
\Big) + S_i \Big] \geq 
\sum_{i=1}^k \Big[{\cal J}(  y_i,  z_i)
- \frac{\min \{ \d_1,\d_2 \}}{8} \Big] \geq
\sum_{i=1}^k \Big[{\cal J}( 0, 0)
- \frac{5\d_3}{8} \Big],
\ee
since $\d_2 < \d_3 $, $| y_i |, |z_i | < \s $ and 
by condition \ref{cond}-$(iii)$. Hence,
since $ {\cal H}_\mu^k (\ov{a},\ov{s})\geq {\cal H}_\mu^k (0, 0, \wtilde{s})$,
by (\ref{Hest})   there exists 
$i_0 \in \{1, \ldots, k \}$ such
that 
\be\label{eq:finalm}
{\cal J}( \ov{a}_2 + y_{i_0}, \ov{a}_3 + z_{i_0}) +
\frac{\d_3}{8}-\frac{|I'_0-I_0|}{k} \ov{a}_3
\geq {\cal J}( 0, 0)- \frac{5 \d_3 }{8}.
\ee
Since $ | a_3 | \leq \rho $ and by (\ref{numberhet}) we get 
$ | I'_0 - I_0 | |\ov{a}_3| / k \leq \d_3/8 $. Hence by 
(\ref{eq:finalm}),
${\cal J}( \ov{a}_2 + y_{i_0}, \ov{a}_3 + z_{i_0})\geq {\cal J}(0,0)-7\d_3/8$.
By condition \ref{cond}-$(iii)$, 
$\ov{a}_3+z_{i_0} \in (-\rho+2\s, \rho-2\s)$ and as a consequence, 
since for all $i$ $|z_i | < \s $, we deduce that  
$\ov{a}_3 + z_i \in (-\rho,\rho)$.
\end{itemize}
We have  proved that the maximum point $(\ov{a}, \ov{s}) \in U$,
which completes the proof of the theorem. 
\end{pf}

As a consequence of the general shadowing theorem 
\ref{step1} and of lemma \ref{step3} we get the following theorem

\begin{theorem}\label{thm:main}
Let $f (\vphi)  =  \sum_{j=1}^n \cos \vphi_j$, $ n \geq 3 $, 
and $ \om_\e $ be a $(\gamma_\e, \tau)$-diophantine vector.  
Assume  $\ep$, $ \mu \e^{-3/2} $ and $ \mu \e^{- 2a -1} $ to be 
sufficiently small. Then, for all   
$ I_0, I_0' $ with $ \om_\e \cdot I_0= \om_\e \cdot I_0' $
and $ (I_0)_1 =  (I_0')_1 $ 
there exists a heteroclinic orbit  connecting the invariant tori
${\cal T}_{I_0}$ and ${\cal T}_{I_0'}$ with a diffusion time  
\be\label{timediffpart}
T_d \leq C \frac{| I_0' - I_0 |}{ \mu \e^{a + (1/2)} } 
\times  \max \Big\{ \frac{1}{\g_\e (\e^{a + (1/2)})^\t }, 
|\ln (\mu ) | \Big\} 
\ee
\end{theorem}

\begin{remark}\label{poliexp}
The number $ k $ of heteroclinic transitions used in the 
proof of theorem \ref{thm:main} is polynomial 
w.r.t $ 1 / \e $ since our shadowing orbit moves along 
the directions $ (I_2, \ldots, I_n) \in {\bf R}^{n-1}$ only
(``directions of large splitting''). On the contrary the shadowing orbit 
connecting tori $ {\cal T}_{I_0} $ and $ {\cal T}_{I_0'} $
with $ ( I_0' )_1 \neq (I_0)_1 $ would shadow 
an exponentially large number of heteroclinic transitions
and the diffusion time would depend also on the contant 
$ \d_2 = 3\pi \mu \e^{-1/2} \exp (-\pi/(2\sqrt{\e}))$ which describes
the exponentially small splitting.
In any case, at each transition, the shadowing orbit approaches   
the homoclinic point only up $ \rho = O( \e^{a + 1/2}) $ and therefore 
the time $ T_s $ spent for each single transition is polynomial w.r.t 
$ 1 / \e $. 
In this way we deduce that the diffusion time $ T_d $ is estimated,
up to inverse powers of $1 / \e $,  by 
an exponential $ T_d = O ( \exp ( \pi/(2\sqrt{\e})) ) $.
Since the (determinant of the) splitting
$ \Delta = O ( \exp ( - \pi / ( 2 \sqrt{\e} ) ) )$ we get that 
$ T_d \approx 1 / \Delta $, while in \cite{BCV} and 
\cite{BB3} the diffusion time is
estimated by
$ T_d \approx 1 / \Delta^p $ for some positive constant $p$.
\end{remark}

\noindent
{\it Massimiliano Berti, S.I.S.S.A., Via Beirut 2-4,
34014, Trieste, Italy, berti@sissa.it}.
\\[2mm]
{\it Philippe Bolle,
D\'epartement de math\'ematiques, Universit\'e
d'Avignon, 33, rue Louis Pasteur, 84000 Avignon, France,
philippe.bolle@univ-avignon.fr}

\end{document}